\documentclass[12pt,reqno]{amsart}
\usepackage{amssymb}

\newtheorem{theorem}{Theorem}
\newtheorem{corollary}{Corollary}

\newtheorem{definition}{Definition}
\newtheorem{example}{Example}

\newtheorem{proposition}{Proposition}
\newtheorem{remark}{Remark}
\numberwithin{equation}{section}
 \numberwithin{remark}{section}
\numberwithin{definition}{section} \numberwithin{lemma}{section}
\numberwithin{proposition}{section}
 \numberwithin{example}{section}
 \numberwithin{corollary}{section}
\numberwithin{theorem}{section}

\begin{document}
\title[]{Rational  self-homotopy  equivalences  and Whitehead exact sequence}
\author{Mahmoud Benkhalifa}
\email{makhalifa@uqu.edu.sa}
\address{Department of Mathematics. Faculty of Applied Sciences. Umm Al-Qura University. Mekka. Saudi Arabia }
\maketitle
\begin{abstract}
For a simply connected   CW-complex $X$, let $\mathcal{E}(X)$ denote
the group of homotopy classes of self-homotopy equivalence of $X$
and let $\mathcal{E}_{\sharp}(X)$ be its subgroup of homotopy
classes which induce the identity on homotopy groups. As we  know,
the quotient group $\frac{\mathcal{E}(X)}{\mathcal{E}_{\sharp}(X)}$
can be identified with a  subgroup  of  $Aut(\pi_{*}(X))$. The aim
of this work is to determine this subgroup for  rational spaces. We
construct the Whitehead exact sequence associated with the minimal
Sullivan model of $X$ which allows us to define the subgroup
$\mathrm{Coh.Aut}(\mathrm{Hom}\big(\pi_{*}(X),\Bbb Q)\big)$ of
self-coherent
 automorphisms of the graded
vector space  $\mathrm{Hom}(\pi_{*}(X),\Bbb Q)$. As a consequence we
establish that
$\frac{\mathcal{E}(X)}{\mathcal{E}_{\sharp}(X)}\cong\mathrm{Coh.Aut}\big(\mathrm{Hom}(\pi_{*}(X),\Bbb
Q)\big)$. In addition, by computing the group
$\mathrm{Coh.Aut}\big(\mathrm{Hom}(\pi_{*}(X),\Bbb Q)\big)$, we give
examples of rational spaces that have few self-homotopy
equivalences.
\end{abstract}

\section{\textbf{Introduction}}
 If $X$ is a simply connected CW-complex, let $\mathcal{E}(X)$
denote the set of homotopy classes of self-homotopy equivalence of
$X$. Equipped with the composition of homotopy classes,
$\mathcal{E}(X)$ is a group. Let $\mathcal{E}_{\sharp}(X)$  denote
the subgroup of homotopy classes which induce the identity on
homotopy groups. Clearly $\mathcal{E}_{\sharp}(X)$ is a normal
subgroup of $\mathcal{E}(X)$. In this paper we study the quotient
group $\frac{\mathcal{E}(X)}{\mathcal{E}_{\sharp}(X)}$ where $X$ is
a rational space. Recall that there exists a homomorphism
$\mathcal{E}(X)\to Aut(\pi_{*}(X))$ whose kernel is
$\mathcal{E}_{\sharp}(X)$, thus we can identify
$\frac{\mathcal{E}(X)}{\mathcal{E}_{\sharp}(X)}$ with a subgroup $G$
of $Aut(\pi_{*}(X))$

The aim of this paper is to determine the subgroup $G$
when $X$ is a  rational space.  Due to the theory elaborated by
Sullivan \cite{H}, the homotopy theory of rational spaces is equivalent to
the homotopy theory of minimal cochain commutative algebras over the
rationals (mccas, for short). Because of this equivalence we can
translate our problem to study the quotient group
$\frac{\mathcal{E}(\Lambda V,\partial)}{\mathcal{E}_{\sharp}(\Lambda
V,\partial)}$, where $(\Lambda V,\partial)$ is  the mcca associated
with $X$ (called the minimal Sullivan model of $X$),
$\mathcal{E}(\Lambda V,\partial)$ denotes the group of self-homotopy
equivalence of $(\Lambda V,\partial)$ and
$\mathcal{E}_{\sharp}(\Lambda V,\partial)$ denotes the subgroup   of
$\mathcal{E}(\Lambda(V),\partial)$ consisting  of the elements
inducing the identity on the indecomposables.

For this purpose  we associate with each mcca $(\Lambda
V,\partial)$ an exact sequence, denoted by $WES(\Lambda V,\partial)$
and
   called the Whitehead exact sequence of
$(\Lambda V,\partial)$. This sequence allows us  to define the
semigroup $\mathrm{Coh.Mor}(V^{*})$ (respect. the group
$\mathrm{Coh.Aut}(V^{*})$)  of self-coherent homomorphisms (respect.
self-coherent automorphisms) of the graded vector space  $V^{*}$ and
to exhibit   a short exact sequence of semigroups (with units): $$
\mathcal{E}_{\sharp}(\Lambda(V),\partial)\rightarrowtail
[(\Lambda(V),\partial),(\Lambda(V),\partial)]\overset{}{\twoheadrightarrow}\mathrm{Coh.Mor}(V^{*})$$
and a short exact sequence of groups:
$$\mathcal{E}_{\sharp}(\Lambda(V),\partial)\rightarrowtail
\mathcal{E}(\Lambda(V),\partial)\overset{}{\twoheadrightarrow}\mathrm{Coh.Aut}(V^{*}),$$
where  $[(\Lambda(V),\partial),(\Lambda(V),\partial)]$  denotes the
semigroup  of homotopy classes
 of cochain morphism from $(\Lambda V ,\partial)$ to itself.

 Because of the homotopy equivalence mentioned above,
 our main result says:\\

 \noindent \textbf{Theorem}. \noindent \emph{If $X$ is a simply connected rational
CW-complex, then:\\
 There exist a short exact sequence of semigroups:
$$\mathcal{E}_{\sharp}(X)\rightarrowtail
[X,X]\overset{}{\twoheadrightarrow}\mathrm{Coh.Mor}\big(\mathrm{Hom}(\pi_{*}(
X),\Bbb Q)\big).$$
 There exist a short exact
sequence of groups:
$$\mathcal{E}_{\sharp}(X)\rightarrowtail
\mathcal{E}(X)\overset{}{\twoheadrightarrow}\mathrm{Coh.Aut}\big(\mathrm{Hom}(\pi_{*}(
X),\Bbb Q)\big).$$} Here $[X,X]$ denotes the semigroup of
self-homotopy maps of  $X$.\\

 In addition and by using techniques of rational homotopy
theory, we compute the groups
$\mathrm{Coh.Aut}\big(\mathrm{Hom}(\pi_{*}( X),\Bbb Q)\big)$ for
certain rational spaces via their minimal Sullivan models. For
instance,  we investigate the following question asked by M. Arkowitz
and G Lupton in \cite{MA}: Which finite groups can be realized as the
group of self-homotopy equivalence of a rational space? We show that
the groups $\underset{2^{n}}{\oplus}\Bbb
 Z_{2}$ ($2^{n}$ copies of $\Bbb
 Z_{2}$) are realizable for $n\leq 10$. At the end of this work we ask the following question: Is it true that   the groups $\underset{2^{n}}{\oplus}\Bbb
 Z_{2}$ are always realizable  for $n\geq 11$?
\section{\textbf{Coherent morphishms}}
\subsection{Whitehead exact sequence of  1-connected mcca}
Let   $(\Lambda V,\partial)$ be a 1-connected mcca. As we have done
in (\cite{ben2}, section 2) \big(respect. (\cite{ben3}, section
2)\big) for 1-connected minimal free cochain algebras (respect. free
chain algebras)  over a P.I.D, we  can define the Whitehead exact
sequence of $(\Lambda V,\partial)$ as follows:

 First define the pair:
$$\Big(\Lambda V^{\leq n+1};\Lambda V^{\leq n-1}\Big)=\frac{(\Lambda V^{\leq n+1},\partial|_{\Lambda V^{\leq n+1}})}{(\Lambda V^{\leq n-1},\partial|_{\Lambda V^{\leq n-1}})}\,\,\,,\,\,\forall n\geq 3$$
where   $\partial|_{\Lambda V^{\leq n}}$ denotes the restriction of
the differential $\partial$ to  $\Lambda V^{\leq n}$.

To each pair $\Big(\Lambda V^{\leq n+1};\Lambda V^{\leq
n-1}\Big)$ corresponds the following short exact sequence of cochain
complexes:
$$(\Lambda V^{\leq n-1}),\partial|_{\Lambda V^{\leq n-1}})\rightarrowtail (\Lambda V^{\leq n+1},\partial|_{\Lambda V^{\leq n+1}})\twoheadrightarrow \Big(\Lambda V^{\leq n+1};\Lambda V^{\leq n-1}\Big)$$
which yields  the following long exact cohomology sequence:

\begin{picture}(300,80)(-80,30)
\put(10,80){$\cdots\rightarrow V^{n}\cong H^{n}\Big(\Lambda V^{\leq
n+1});\Lambda V^{\leq n-1})\Big)\overset{b^{n}}{\longrightarrow}
H^{n+1}(\Lambda V^{\leq n-1})) $} \put(263,76){$\vector(0,-1){24}$}
 \put(-92,40){$\cdots \leftarrow H^{n+1}(\Lambda V^{\leq
n+1}))\leftarrow V^{n+1}\cong   H^{n+1}\Big(\Lambda V^{\leq
n+1});\Lambda V^{\leq n-1})\Big) \overset{j^{n+1}}{\longleftarrow}
H^{n+1}(\Lambda V^{\leq n+1}))$}
\end{picture}

 Consequently,  if we combine the two long exact cohomology
sequences associated with the two pairs $\Big(\Lambda V^{\leq
n+1};\Lambda V^{\leq n-1}\Big)$ and $\Big(\Lambda V^{\leq
n+2};\Lambda V^{\leq n}\Big)$ respectively, we get the following long
exact  sequence:
\begin{eqnarray}
\label{23} \rightarrow V^{n}\overset{b^{n}}{\rightarrow}
H^{n+1}(\Lambda V^{\leq n-1}) \rightarrow H^{n+1}(\Lambda V^{\leq
n+1})\overset{j^{n+1}}{\rightarrow}   V^{n+1}
\overset{b^{n+1}}{\rightarrow} H^{n+2}(\Lambda V^{\leq
n})\rightarrow \,\,\,\,\,
\end{eqnarray}
where the homomorphisms $b ^{n}$ and $j^{n+1}$ are defined as
follows:
\begin{equation}
b^{n}(v_{n})=[\partial^{n}(v_{n})]\hspace{2cm}
j^{n+1}([v_{n+1}+q_{n+1}])=v_{n+1}. \label{106}
\end{equation}
Here $[\partial^{n}(v_{n})]$ and $[v_{n+1}+q_{n+1}]$ denote
respectively the cohomology classes of $\partial^{n}(v_{n})\in
(\Lambda V^{\leq n-1})^{n+1}$ and  $v_{n+1}+q_{n+1}\in (\Lambda
V^{\leq n+1})^{n+1}$.

 Since it is well-known that $H^{n+1}(\Lambda V^{\leq
n+1})\cong H^{n+1}(\Lambda V)$, then from  (\ref{23}) we get  the
following long  sequence:
\begin{equation*}
\cdots \rightarrow V^{n}\overset{b^{n}}{\longrightarrow }%
H^{n+1}(\Lambda V^{\leq n-1})\longrightarrow H^{n+1}(\Lambda
V)\longrightarrow
V^{n+1}%
\overset{b^{n+1}}{\longrightarrow }\cdots
\end{equation*}
 called  the Whitehead exact sequence of $(\Lambda
V,\partial)$.

This  sequence is natural with respect to  cochain morphisms. That
is, if $\alpha:(\Lambda(V),\partial )\rightarrow (\Lambda(W),\delta
)$ is a cochain morphism, then $\alpha$ induces the following
commutative diagram:

\begin{picture}(300,90)(25,15)
\put(65,80){$\cdots \rightarrow V^{n}\overset{b^{n}}{\longrightarrow
}H^{n+1}(\Lambda V^{\leq n-1} )\longrightarrow H^{n+1}(\Lambda V
)\longrightarrow V^{n+1}\overset{b^{n+1}}{\longrightarrow }\cdots$}
\put(257,76){$\vector(0,-1){44}$} \put(147,76){$\vector(0,-1){44}$}
 \put(62,20){$\cdots \rightarrow
W^{n}\overset{b'^{n}}{\longrightarrow }H^{n+1}(\Lambda W^{\leq n-1}
)\longrightarrow H^{n+1}(\Lambda W )\longrightarrow
W^{n+1}\overset{b'^{n+1}}{\longrightarrow }\cdots$}
\put(100,76){$\vector(0,-1){44}$} \put(88,50){\scriptsize
$\widetilde{\alpha}^{n}$} \put(332,50){\scriptsize
$\widetilde{\alpha}^{n+1}$} \put(330,76){$\vector(0,-1){44}$}
\put(148,50){\scriptsize $H^{n+1}( \alpha_{(n-1)})$}
\put(257,50){\scriptsize  $H^{n+1}(\alpha)$} \put(10,50){
\scriptsize $(1)$}
\end{picture}

\noindent where $\widetilde{\alpha}:V^{*}\to W^{*}$ is the graded
homomorphism induced by $\alpha$ on the indecomposables and where
$\alpha_{(n-1)}:(\Lambda V^{\leq n-1} ,\partial)\to (\Lambda W^{\leq
n-1} ,\delta)$ is the restriction of $\alpha$.

\subsection{Coherent morphisms between  Whitehead exact sequences}
 Let $(\Lambda V ,\partial)$, $ (\Lambda W,\delta)$ be two
1-connected mccas and let  $\xi:V^{*}\to W^{*}$ be a given graded
linear application. For every $n\geq 2$, let $\{\xi^{\leq n}\}$
denote
 the set of all cochain morphisms
from  $(\Lambda V^{\leq n} ,\partial)$ to $(\Lambda W^{\leq
n} ,\delta)$ inducing $\xi^{\leq n}$ on the indecomposables.\\
\begin{definition}
\label{d1} \textnormal{Let $(\Lambda V ,\partial)$ and $(\Lambda
W,\delta)$ be two 1-connected mccas. A graded linear map
$\xi^{*}:V^{*}\to W^{*}$ is called  a coherent morphism   if the
following holds}:

\noindent \textnormal{For every $n\geq 2$, if the set $\{\xi^{\leq
n}\}$ is not empty, then it contains $\alpha_{(n)}$ making  the
following diagram commute}:

\begin{picture}(300,90)(-30,30)
\put(60,100){$V^{n+1}\hspace{1mm}\vector(1,0){150}\hspace{1mm}W^{n+1}$}
 \put(68,76){\scriptsize $b^{n+1}$} \put(248,76){\scriptsize $b'^{n+1}$}
\put(66,96){$\vector(0,-1){37}$} \put(245,96){$\vector(0,-1){37}$}
\put(155,103){\scriptsize $\xi^{n+1}$} \put(145,52){\scriptsize
$H^{n+2}(\alpha_{(n)})$} \put(35,48){$H^{n+2}(\Lambda V^{\leq n} )
\hspace{1mm}\vector(1,0){130}\hspace{1mm}H^{n+2}(\Lambda W^{\leq
n})\hspace{1mm}$} \put(-29,76){\scriptsize $(2)$}
\end{picture}
\end{definition}
\begin{example}
\label{e1} \textnormal{If $\alpha:(\Lambda(V),\partial )\rightarrow
(\Lambda(W),\delta )$ is a cochain algebra morphism between two 1-connected
mccas, then, according to diagram $(1)$, the  graded linear map
$\widetilde{\alpha}:V^{*}\to W^{*}$  is a coherent morphism}.
\end{example}
\begin{example}
\label{e5} \textnormal{It easy to see that if $(\Lambda(V),\partial)
$ is a 1-connected mcca, then $Id_{V^{*}}$ is a coherent morphism.
Observe that in this case the
 set $\{Id_{V^{*}}\}$ of cochain morphisms from $(\Lambda(V),\partial)
 $ to it self inducing $Id_{V^{*}}$ on the indecomposables is always not empty since it contains $Id_{(\Lambda(V),\partial)}$}.
\end{example}

Now let $(\Lambda V ,\partial)$, $ (\Lambda W,\delta)$ be two
1-connected mccas and let  $\mathrm{Coh.Mor}(V^{*},W^{*})$ denote
the set of all the coherent automorphisms from $V^{*}$  to $W^{*}$.
Example \ref{e1} allows us to define a map $\Phi:[(\Lambda V
,\partial),(\Lambda W,\delta)]\to \mathrm{Coh.Mor}(V^{*},W^{*})$ by
setting $\Phi([\alpha])=\widetilde{\alpha}$. Here $[(\Lambda V
,\partial),(\Lambda W,\delta)]$ denote the set of homotopy classes
 from $(\Lambda V ,\partial)$ to $(\Lambda W ,\delta)$. Recall that there is a reasonable concept of ``homotopy''  among cochain
 morphisms  (see for example \cite{H} for details),
 analogous in many respects to the topological notion of homotopy.
 \begin{remark}
\label{r1} \textnormal{It is well-known (\cite{H} proposition 12.8)
that if two cochain morphisms $\alpha,\alpha':(\Lambda V
,\partial)\to(\Lambda W,\delta)$ are homotopic, then they induce the
same graded linear maps  on the indecomposables i.e,
$\widetilde{\alpha}=\widetilde{\alpha'}$. So the map $\Phi$ is
well-defined}.\\
\end{remark}

\begin{proposition}
\label{p1} The map $\Phi$ is surjective.\\
\end{proposition}

\noindent \textbf{Proof.} Let $\xi\in\mathrm{Coh.Mor}(V^{*},W^{*})$.
Assume, by induction, that we have constructed  a cochain morphism
$\theta_{(n)}: (\Lambda V^{\leq n} ,\partial)\to (\Lambda W^{\leq
n},\delta)$ such that $\Phi([\theta_{(n)}])=\xi^{\leq n}$. This
implies  that the set $\{\xi^{\leq n}\}$ is not empty.  Therefore,
by definition \ref{d1}, this set contains an element $\alpha_{(n)}$
making the diagram (2) commutes. Now choose
$(v_{\sigma})_{\sigma\in\Sigma}$ as a basis of $V^{n+1}$. Recall
that, in this diagram, we have:
\begin{eqnarray}
\label{117} H^{n+2}(\alpha_{(n)})\circ
b^{n+1}(v_{\sigma})&=&\alpha_{(n)}\circ\partial^{n+1}(v_{\sigma})+\mathrm{Im}\,\delta^{n+1}\nonumber\\
b'^{n+1}\circ
\xi^{n+1}(v_{\sigma})&=&\delta^{n+1}\circ\xi^{n+1}(v_{\sigma})+\mathrm{Im}\,\delta^{n+1}\label{103}
\end{eqnarray}
where $\delta^{n+1}:(\Lambda W^{\leq n})^{n+1} \to (\Lambda W^{\leq
n})^{n+2} $. Since the diagram $(2)$ commutes,  the element
$(\alpha_{(n)}\circ\partial^{n+1}-\delta^{n+1}\circ\xi^{n+1})(v_{\sigma})\in\mathrm{Im}
\,\delta^{n+1}$. As a consequence there exists $u_{\sigma}\in
(\Lambda W^{\leq n})^{n+1}$ such that:
\begin{eqnarray}
(\alpha_{(n)}\circ\partial^{n+1}-\delta^{n+1}\circ\xi^{n+1})(v_{\sigma})=\delta^{n+1}(u_{\sigma})\label{123}.
\end{eqnarray}
Thus we define $\theta_{(n+1)}:(\Lambda V^{\leq n+1}
,\partial)\rightarrow(\Lambda W^{\leq n+1},\delta)$ by setting:
\begin{eqnarray}
\theta_{(n+1)}(v_{\sigma})=\xi^{n+1}(v_{\sigma})+u_{\sigma}\,\,\,,\,\,
v_{\sigma}\in V^{n+1}\,\,\,\,\,\text{and}\,\,\,\,\,
\theta=\alpha_{(n)}\text{ on }V^{i}\,\,\,,\,\,\forall i\leq
n\label{113}
\end{eqnarray}
As $\partial^{n+1}(v_{\sigma})\in (\Lambda V^{\leq n})^{n+2}$ then,
by (\ref{123}), we get:
\begin{eqnarray*}
\delta^{n+1}\circ\theta_{(n+1)}(v_{\sigma})=\delta^{n+1}(\xi^{n+1}(v_{\sigma}))+\delta^{n+1}(u_{\sigma})=\alpha_{(n)}\circ\partial^{n+1}(v_{\sigma})=\alpha_{(n)}\circ\partial^{n+1}(v_{\sigma}).
\end{eqnarray*}
So  $\theta_{(n+1)}$ is a  cochain morphism.   Now due to the fact
that $u_{\sigma}\in (\Lambda W^{\leq n})^{n+1}$,   the homomorphism
$\widetilde{\theta}_{(n+1)}^{n+1}:V^{n+1}\to W^{n+1}$ coincides with
$\xi^{n+1}$. This implies that $\Phi([\theta_{(n+1)}])=\xi^{\leq
n+1}$ and the set $\{\xi^{\leq n+1}\}$ is not empty, completing the
induction step. Finally the iteration of this process yields a
cochain morphism $\theta:(\Lambda V ,\partial)\to(\Lambda W,\delta)$
satisfying
$\widetilde{\theta}=\xi$\hspace{1cm}$\square$\\
\begin{remark}
\label{r3} \textnormal{If we assume that $(\Lambda W^{\leq
n})^{n+1}=0$ for all $n\leq k$,  then  the cochain morphism
$\theta_{(n+1)}$ given in (\ref{113}) will satisfy
$\theta_{(n+1)}=\xi^{\leq n+1}$  and the set $\{\xi^{\leq n+1}\}$
contains just one element for all $n\leq k$}.\\
\end{remark}
\begin{remark}
\label{r2} \textnormal{It  is well-known (see \cite{H}) that  any
cochain morphism  between two 1-connected mccas inducing a graded
linear isomorphism on the indecomposables is an isomorphism.
Consequently if the coherent morphism $\xi$ is an isomorphism, then
the  cochain morphism $\alpha:(\Lambda V ,\partial)\to(\Lambda
W,\delta)$ constructed  in the proof of proposition \ref{p1} is such
that $\alpha_{( n)}:(\Lambda V^{\leq n} ,\partial)\to(\Lambda
W^{\leq n},\delta)$ is a cochain isomorphism
for every $n\geq 2$}.
\end{remark}
Now  let us denote by $\mathcal{E}(\Lambda(V),\partial)$ the group
of the self-homotopy equivalences of $(\Lambda(V),\partial)$ and by
$\mathcal{E}_{\sharp}(\Lambda(V),\partial)$ the subgroup of
$\mathcal{E}(\Lambda(V),\partial)$ consisting  of the elements
inducing the identity on the  indecomposables. Also let
$\mathrm{Coh.Aut}(V^{*})$ denote the set of the self-coherent
automorphisms of $V^{*}$.
\begin{proposition}
\label{p2} $\mathrm{Coh.Aut}(V^{*})$ is a  subgroup of the  group
$Aut(V^{*})$.
\end{proposition}

\noindent \textbf{Proof.}  Let  $\xi,\xi'\in
\mathrm{Coh.Aut}(V^{*})$. By definition \ref{d1} to prove that
$\xi'\circ\xi\in \mathrm{Coh.Aut}(V^{*})$, we must show that, for
every $n\geq 2$, if the set $\{(\xi'\circ\xi)^{\leq n}\}$ is not
empty, then it contains an element $\lambda_{(n)}$ making making the
following diagram commutes:

\begin{picture}(300,90)(-30,30)
\put(60,100){$V^{n+1}\hspace{1mm}\vector(1,0){150}\hspace{1mm}V^{n+1}$}
 \put(68,76){\scriptsize $b^{n+1}$} \put(248,76){\scriptsize $b^{n+1}$}
\put(66,96){$\vector(0,-1){37}$} \put(245,96){$\vector(0,-1){37}$}
\put(155,103){\scriptsize $\xi'^{n+1}\circ\xi^{n+1}$}
\put(145,52){\scriptsize $H^{n+2}(\lambda_{(n)})$}
\put(35,48){$H^{n+2}(\Lambda V^{\leq n} )
\hspace{1mm}\vector(1,0){130}\hspace{1mm}H^{n+2}(\Lambda V^{\leq
n})\hspace{1mm}$}
\end{picture}

Recall that $\{(\xi'\circ\xi)^{\leq n}\}$  is the set of
all cochain morphisms from  $(\Lambda V^{\leq n} ,\partial)$ to
itself  inducing $(\xi'\circ\xi)^{\leq n}$ on the indecomposables.

 Indeed, if  $\xi,\xi'\in \mathrm{Coh.Aut}(V^{*})$, then,
according the proof of  proposition \ref{p1}, there  exist two
cochain isomorphisms $\alpha,\alpha':(\Lambda(V),\partial)\to
(\Lambda(V),\partial)$ such that $\widetilde{\alpha}=\xi$,
$\widetilde{\alpha}'=\xi'$ and satisfying:
\begin{equation*}
  b^{n+1}\circ \xi'^{n+1}=H^{n+2}(\alpha'_{(n)})\circ b^{n+1}\hspace{1cm} b^{n+1}\circ \xi^{n+1}=H^{n+2}(\alpha_{(n)})\circ
  b^{n+1}\,\,\,\,\,,\,\,\,\,\,\forall n\geq 2
\end{equation*}
So, for all $n\geq 2$,  the set $\{(\xi'\circ\xi)^{\leq n}\}$
contains $\lambda_{(n)}=\alpha'_{(n)}\circ\alpha_{(n)}$  and an easy
computation shows that:
\begin{eqnarray}\label{2}
  b^{n+1}\circ \xi'^{n+1}\circ \xi^{n+1}&\hspace{-2mm}=&\hspace{-2mm}H^{n+2}(\alpha'_{(n)})\circ b^{n+1}\circ
  \xi^{n+1}=H^{n+2}(\alpha'_{(n)})\circ H^{n+2}(\alpha'_{(n)})\circ b^{n+1}\nonumber\\
&\hspace{-2mm}=&\hspace{-2mm}H^{n+2}(\alpha'_{(n)}\circ\alpha'_{(n)})\circ
b^{n+1}\,=\,H^{n+2}(\lambda_{(n)})\circ b^{n+1}
\end{eqnarray}
\noindent which implies that    $\xi'\circ \xi$ is a coherent
automorphism.

Now  let  $\xi\in \mathrm{Coh.Aut}(V^{*})$. By proposition
\ref{p1} and remark \ref{r2} we get  a   cochain isomorphism
$\alpha:(\Lambda V,\partial)\to (\Lambda V ,\partial)$
  satisfying  $\widetilde{\alpha}=\xi$ and such that   $\alpha_{(n)}$ is  a cochain isomorphism  for all $n\geq
  2$. Consequently there exists $\alpha'_{(n)}$ such that
$\alpha_{(n)}\circ\alpha'_{(n)}=\alpha'_{(n)}\circ\alpha_{(n)}=Id_{(\Lambda
V^{\leq n},\partial)}$ which implies that
$\widetilde{\alpha'}_{(n)}=(\xi^{-1})^{n}$ for all $n\geq
  2$. So the set
$\{(\xi^{-1})^{\leq n}\}$ contains
  $\alpha'_{(n)}$. Moreover as  $\xi\in \mathrm{Coh.Aut}(V^{*})$ it  satisfies $b^{n+1}\circ \xi^{n+1}=H^{n+2}(\alpha_{(n)})\circ
b^{n+1}$ which implies that $H^{n+2}(\alpha'_{(n)})\circ
b^{n+1}=b^{n+1}\circ (\xi^{n+1})^{-1}$. Hence   $\xi^{-1}\in
\mathrm{Coh.Aut}(V^{*})$\hspace{1cm}$\square$\\

\begin{theorem}
\label{t2}Let  $(\Lambda(V),\partial)$ be a  1-connected mcca. There
exists a short exact sequence of groups:
\begin{equation}
\mathcal{E}_{\sharp}(\Lambda(V),\partial)\rightarrowtail
\mathcal{E}(\Lambda(V),\partial)\overset{\Phi}{\twoheadrightarrow}\mathrm{Coh.Aut}(V^{*})\label{105}\\
\end{equation}
\end{theorem}

\noindent \textbf{Proof.} First we have $\Phi(\{\alpha\circ
\alpha'\})=\widetilde{\alpha\circ \alpha'}=\widetilde{\alpha}\circ
\widetilde{\alpha}'=\Phi(\{\alpha\})\circ \Phi(\{\alpha'\})$. Next
the surjection of $\Phi$  is assured by  proposition  \ref{p1} and
finally it is clear that $\ker \Phi=\mathcal{E}_{\sharp}(\Lambda V
,\partial)$ \hspace{1cm}$\square$\\

  The set $[(\Lambda(V),\partial),(\Lambda(V),\partial)]$ of self-homotopy
classes of a 1-connected mcca $(\Lambda V,\partial)$, equipped with
the composition of maps, is a semigroup with unit. So let
$\mathrm{Coh.Mor}(V^{*})$ denote the set of the self-coherent
morphisms of $V^{*}$. From proposition \ref{p1} we deduce that
$\mathrm{Coh.Mor}(V^{*})$ is a semigroup with unit and the map
$\Phi$ is a homomorphism of semigroups. Hence
theorem \ref{t2} implies that:\\
\begin{corollary}
\label{c2} Let  $(\Lambda(V),\partial)$ be a  1-connected mcca.
There exists a short exact sequence of semigroups:
\begin{equation}
\mathcal{E}_{\sharp}(\Lambda(V),\partial)\rightarrowtail
[(\Lambda(V),\partial),(\Lambda(V),\partial)]\overset{\Phi}{\twoheadrightarrow}\mathrm{Coh.Mor}(V^{*}).\label{55}
\end{equation}
\end{corollary}
Because of the equivalence between  the homotopy theory of rational
spaces and  the homotopy theory of mccas, we can construe  the above
results as follows. Let $X$ be a simply connected rational
CW-complex of finite type. By  the properties of  the Sullivan minimal
 model $(\Lambda(V),\partial)$  of $X$, we can identify $\mathcal{E}(X)$  with  $\mathcal{E}(\Lambda(V),\partial)$  and $\mathcal{E}_{\sharp}(X)$ with $\mathcal{E}_{\sharp}(\Lambda(V),\partial)$.
 Moreover $WES(\Lambda(V),\partial)$ can be written as follows:
\begin{equation*}
 \cdots \rightarrow \mathrm{Hom}(\pi_{n}( X),\Bbb Q)\overset{b_{X}^{n-1}}{\rightarrow }%
\Gamma_{X} ^{n}\rightarrow H^{n}(X,\Bbb Q)\rightarrow
\mathrm{Hom}(\pi_{n}( X),\Bbb Q)%
\overset{b_{X}^{n}}{\rightarrow }\cdots
\end{equation*}
where $\Gamma ^{n}_{X}=H^{n}(\Lambda(V^{\leq n-2}))$. We call this
sequence  the Whitehead exact sequence of $X$ and we denote it by
$WES(X)$. Clearly  this sequence is an invariant of homotopy.

 As a  consequences of  theorem \ref{t2} and corollary
\ref{c2} we establish the following result:\\
\begin{corollary}
\label{108} There exist a short exact sequence of groups:
\begin{equation}
\mathcal{E}_{\sharp}(X)\rightarrowtail
\mathcal{E}(X)\overset{\Psi}{\twoheadrightarrow}\mathrm{Coh.Aut}\big(\mathrm{Hom}(\pi_{*}(
X),\Bbb Q)\big).\label{121}
\end{equation}
There exist a short exact sequence of semigroups:
\begin{equation}
\mathcal{E}_{\sharp}(X)\rightarrowtail
[X,X]\overset{\Psi}{\twoheadrightarrow}\mathrm{Coh.Mor}\big(\mathrm{Hom}(\pi_{*}(
X),\Bbb Q)\big)\label{121}
\end{equation}
\end{corollary}
Here $[X,X]$ denotes the semigroup of the self-homotopy classes of
$X$ and  the linear map $\Psi(\{\alpha\}):\mathrm{Hom}\big(\pi_{*}(
X),\Bbb Q\big)\to \mathrm{Hom}\big(\pi_{*}( X),\Bbb Q\big)$ is
defined as follows:
$$\Psi(\{\alpha\})(\eta)=\eta\circ \pi_{*}(\alpha)\,\,\,,\,\,\,\forall \eta\in \mathrm{Hom}\big(\pi_{*}( X),\Bbb
Q\big)$$

  Recall that  $\mathrm{Coh.Aut}\big(\mathrm{Hom}(\pi_{*}(
 X),\Bbb Q)\big)$ (respect. $\mathrm{Coh.Mor}\big(\mathrm{Hom}(\pi_{*}(
X),\Bbb Q)$\big) denotes the subgroup of the self-coherent
automorphisms (respect. self-coherent morphisms)
 of  $\mathrm{Hom}(\pi_{*}(
X),\Bbb Q)$.\\
\begin{definition}
\label{d2} \textnormal{Let  $(\Lambda V ,\partial)$ and $ (\Lambda
W,\delta)$ be two 1-connected mccas. We say that $WES(\Lambda V
,\partial)$ and $ WES(\Lambda W,\delta)$ are coherently isomorphic
if the set $\mathrm{Coh.Iso}(V^{*},W^{*})$ of the coherent isomorphisms from $V^{*}$ to $W^{*}$  is not empty}.\\
\end{definition}
\begin{corollary}
\label{c3}Two simply connected rational CW-complexes of finite type
$X$ and $Y$  are homotopy equivalent if and only their  $WES(X)$  and $WES(Y)$
are coherently
isomorphic.\\
\end{corollary}
\noindent \textbf{Proof.} Let $(\Lambda V,\partial)$ (respect.
$(\Lambda W,\delta)$)  the Sullivan minimal
 model of $X$ (respect. of $Y$). First recall that $WES(X)=WES(\Lambda V,\partial)$  and $WES(Y)=WES(\Lambda W,\delta)$. Now if  $WES(X)$  and $WES(Y)$ are coherently
isomorphic, then there exists a coherent isomorphism $\xi: V^{*}\to
 W^{*}$. Now proposition \ref{p1} yields a cochain
morphism $\alpha:(\Lambda V,\partial)\to (\Lambda W,\delta)$ such
that the map $\widetilde{\alpha}$, induced by $\alpha$ on the
indecomposables, satisfies $\widetilde{\alpha}=\xi$. So
$\widetilde{\alpha}$ is an isomorphism which means that the models
$(\Lambda V,\partial)$ and $(\Lambda W,\delta)$ are isomorphic.
Hence,  by the properties of the Sullivan minimal
 model, we conclude that
$X$ and $Y$ are homotopy equivalent \hspace{1cm}$\square$\\
\section{\textbf{Examples}}
In this section we give some examples showing how the group
$\mathrm{Coh.Aut} \big(\mathrm{Hom}(\pi_{*}(
 X),\Bbb Q)\big)$ can be used to compute the group $\mathcal{E}(X)$ when $X$ is a simply connected rational
 CW-complex. First let us consider the following example which has already
treated  in (\cite{MA}, example 5.3), where the authors have used
another technique, which is radically different from our approach,
to
determine $\mathcal{E}(X)$.\\

\begin{example}
\label{e3}  \textnormal{Let $\Lambda V=\Lambda
(x_{1},x_{2},y_{1},y_{2},y_{3},z)$ with $|x_{1}|=10$, $|x_{2}|=12$,
$|y_{1}|=41$, $|y_{2}|=43$, $|y_{3}|=45$ and $|z|=119$. The
differential is as follows}:

$
\begin{array}{ccccccl}
  \partial(x_{1})=0 &&& \partial(y_{1})=x_{1}^{3}x_{2} &&&\partial(y_{3})=x_{1}x_{2}^{3} \\
  \partial(x_{2})=0 &&&\partial(y_{2})=x_{1}^{2}x_{2}^{2}  &&&
 \,\, \partial(z)=y_{1}y_{2}x_{2}^{3}-y_{1}y_{3}x_{1}x_{2}^{2}+y_{2}y_{3}x_{1}^{2}x_{2}+x_{1}^{12}+x_{2}^{10}
\end{array}
$

\noindent \textnormal{ An easy computation shows that:
$$ H^{46}(\Lambda V^{\leq 44}
)= \Bbb Q\{x_{1}x_{2}^{3}\}
 \hspace{1cm} H^{44}(\Lambda V^{\leq 42}) = \Bbb Q\{x_{1}^{2}x_{2}^{2}
  \}\hspace{1cm}
  H^{42}(\Lambda V^{\leq 40}) = \Bbb Q\{x_{1}^{3}x_{2}\}$$
  $$H^{120}(\Lambda V^{\leq 118})=  \Bbb
Q\{y_{1}y_{2}x_{2}^{3}-y_{1}y_{3}x_{1}x_{2}^{2}+y_{2}y_{3}x_{1}^{2}x_{2},x_{1}^{12},x_{2}^{10}\}$$
First  any linear map $\xi^{i}:V^{i}\to V^{i}$, where
$i=10,12,41,43,45,119$, is  multiplication with a rational number,
so write $\xi^{i}=p_{i}$. Hence in this case any element of
$\mathrm{Coh.Mor}(V^{*})$ can be identified with
$(p_{10},p_{12},p_{41},p_{43},p_{45},p_{119})\in \Bbb Q^{6},$
therefore  $\mathrm{Coh.Mor}(V^{*})$ can be regarded as a semigroup
of $(\Bbb Q^{6},\times)$}.

\noindent \textnormal{Now define the cochain algebra morphism
$\alpha_{(40)}:\Lambda V^{\leq 40}\to \Lambda V^{\leq 40}$ by
 $\alpha_{(40)}(x_{1})=\xi^{10}(x_{1})$  and $\alpha_{(40)}(x_{2})=\xi^{12}(x_{2})$. So the set  of the cochain morphisms from
$(\Lambda V^{\leq 40},\partial)$ to itself inducing
$\xi^{10},\xi^{12}$ on the indecomposables is not empty. To be a
coherent morphism the linear map $\xi^{41}:V^{41}\to V^{41}$ must
satisfy, according to the diagram (2), the relation:
\begin{equation}\label{44}
     b^{41}\circ
\xi^{41}=H^{42}(\alpha_{(40)})\circ b^{41}
\end{equation}
where the linear map  $b^{41}:V^{41}=\Bbb Q\{y_{1}\}\to
H^{42}(\Lambda V^{\leq 40}) = \Bbb Q\{x_{1}^{3}x_{2}\}$ can  be
regarded  as multiplication with a rational. Now as
$H^{42}(\alpha_{(40)})$ is identified with  multiplication by
$p_{10}^{3}p_{12}$, the relation (\ref{44}) implies the equation
$p_{41}=p_{10}^{3}p_{12}$. By going back to the proof of proposition
\ref{p1},  this equation allows us to extend
$\alpha_{(40)}$ to  a cochain morphism $\alpha_{(41)}:\Lambda
V^{\leq 41}\to \Lambda V^{\leq 41}$. Because   $(\Lambda V^{\leq
40})^{ 41}=0$, the element $u_{\sigma}\in (\Lambda V^{\leq
40})^{41}$, given in (\ref{123}), is zero. Consequently   we have
$\alpha_{(41)}(y_{1})=\xi^{41}(y_{1})$.}

 \textnormal{By using a  similar argument in degree
$43,45,119$ we get the following equations:
$$p_{43}=p_{10}^{2}p_{12}^{2}\,\,\,,\,\,\,p_{45}=p_{10}p_{12}^{3}\,\,\,,\,\,\,p_{119}=p_{41}p_{43}p_{10}^{2}p_{12}^{3}=p_{10}^{12}=p_{12}^{10}$$
which have  3 solutions $(0,0,0,0,0,0)$, $(1,1,1,1,1,1)$,
$(1,-1,-1,1,-1,1)$. So we get 3 coherent homomorphisms.}

 \textnormal{ Now by corollary \ref{t2} we have
$\Phi^{-1}(Id_{V_{*}})=\xi_{\sharp}(\Lambda V,\partial)$.  Due to
the fact that $(\Lambda V^{\leq n})^{ n+1}=0$ for
$n=10,12,41,43,45,119$ and by using remark (\ref{r3},   we deduce
that $\Phi^{-1}(Id_{V_{*}})=\{Id_{(\Lambda V,\partial)}\}$ which
implies that $\Phi$, given in the short exact sequence (\ref{55}),
is an isomorphism of semigroups. Hence $[(\Lambda
V,\partial),(\Lambda V,\partial)]$ has 3 elements and then
$\mathcal{E}(\Lambda V,\partial)$ has 2 elements corresponding to
the coherent automorphisms $(1,1,1,1,1,1)$, $(1,-1,-1,1,-1,1)$}.
\end{example}
\begin{example}
\label{e4} \textnormal{Let $(\Lambda W,\delta)$ be the mcca
obtained from the graded algebra $(\Lambda V,\partial)$, given in
example \ref{e3}, by adding  a new generator $|x_{0}|$, with
 $|x_{0}|=2$ and where the differential is as follows, $\delta(x_{0})=0$, $\delta(z)=y_{1}y_{2}x_{2}^{3}-y_{1}y_{3}x_{1}x_{2}^{2}+y_{2}y_{3}x_{1}^{2}x_{2}+x_{1}^{12}+x_{2}^{10}+x_{0}^{60}$ and $\delta=\partial$ on the other generators.
 In this case a simple computation shows that:
$$ \mathrm{Im}\,b^{46}= \Bbb Q\{x_{1}x_{2}^{3}\}\,\, ,\,\,  \mathrm{Im}\,b^{44} = \Bbb Q\{x_{1}^{2}x_{2}^{2}
  \}\,\, ,\,\,
 \mathrm{Im}\,b^{42}= \Bbb Q\{x_{1}^{3}x_{2}\}$$
  $$\mathrm{Im}\,b^{120}=  \Bbb
Q\{y_{1}y_{2}x_{2}^{3}-y_{1}y_{3}x_{1}x_{2}^{2}+y_{2}y_{3}x_{1}^{2}x_{2},x_{1}^{12},x_{2}^{10},x_{0}^{60}\}$$
Write $\xi^{2}:W^{2}=\Bbb Q\{x_{0}\}\to W^{2}=\Bbb Q\{x_{0}\}$ as
$\xi^{2}(x_{0})=p_{2}x_{0}$ with $p_{2}\in \Bbb Q$. By similar
arguments as in example \ref{e3} we get the following equations:
$$p_{41}=p_{10}^{3}p_{12}\,\,\,,\,\,\,p_{43}=p_{10}^{2}p_{12}^{2}\,\,\,,\,\,\,p_{45}=p_{10}p_{12}^{3}\,\,\,,\,\,\,p_{119}=p_{41}p_{43}p_{10}^{2}p_{12}^{3}=p_{10}^{12}=p_{12}^{10}=p_{2}^{60}$$
which give  5 coherent homomorphisms:
$$(0,0,0,0,0,0,0),(1,1,1,1,1,1,1),(1,1,-1,-1,1,-1,1)$$
$$(-1,1,1,1,1,1,1),(-1,1,-1,-1,1,-1,1).$$
As in the example \ref{e3} we have
$\Phi^{-1}(Id_{W_{*}})=\{Id_{(\Lambda W,\delta)}\}$, so $[(\Lambda
W,\delta),(\Lambda W,\delta)]$ has 5 elements and then
$\mathcal{E}(\Lambda W,\delta)$ has 4 elements corresponding to the
coherent automorphisms:
$$(1,1,1,1,1,1,1),
(1,-1,-1,1,-1,1),(-1,1,1,1,1,1,1),(-1,1,-1,-1,1,-1,1).$$ As the last
three elements are of order 2
 we conclude that  $\mathcal{E}(\Lambda W,\delta)\cong \Bbb Z_{2}\oplus\Bbb
 Z_{2}$}.
\end{example}
\begin{remark}
\label{r5}\textnormal{ In \cite{MA} M. Arkowitz and G Lupton ended
their work by the following question: Which finite groups can be
realized as the group of self-homotopy equivalence of a rational
space? Examples \ref{e3} and \ref{e4} show that the groups $\Bbb
Z_{2}$ and $\Bbb Z_{2}\oplus\Bbb
 Z_{2}$ are realizable}.

 \noindent \textnormal{Now let $(\Lambda U_{1},\partial_{1})$ be the  mcca obtained from $(\Lambda W,\delta)$, given in example \ref{e4}, by adding a new generator $x_{3}$  with
 $|x_{3}|=3$ and where the differential is $\partial_{1}(x_{3})=0$, $\partial_{1}(z)=\delta(z)+x_{3}^{40}$ and $\partial_{1}=\delta$ on the other generators.
If we write $\xi^{3}:U_{1}^{3}=\Bbb Q\{x_{3}\}\to U_{1}^{3}=\Bbb
Q\{x_{3}\}$ as $\xi^{3}(x_{3})=p_{3}x_{3}$ with $p_{3}\in \Bbb Q$,
then we will get the following  equations:
\begin{equation}\label{1}
p_{41}=p_{10}^{3}p_{12}\,\,\,,\,\,\,p_{43}=p_{10}^{2}p_{12}^{2}\,\,\,,\,\,\,p_{45}=p_{10}p_{12}^{3}\,\,\,,\,\,\,p_{119}=p_{41}p_{43}p_{10}^{2}p_{12}^{3}=p_{10}^{12}=p_{12}^{10}=p_{2}^{60}=p_{3}^{40}
\end{equation}
which have  the following nontrivial solutions:
$$p_{2}=p_{3}=p_{12}=p_{41}=p_{45}=\pm
1\,\,\,\,\,\,\,\,,\,\,\,\,\,\,\,\,p_{10}=p_{43}=p_{119}=1.$$ which
give  8 coherent automorphisms (seven of them are of order 2) which
are:
$$(1,1,1,1,1,1,1,1),(1,1,1,1,1,1,1,-1),(1,1,-1,-1,1,-1,1,1),(1,1,-1,-1,1,-1,1,-1)$$
$$(-1,1,1,1,1,1,1,1),(-1,1,1,1,1,1,1,1,-1),(-1,1,-1,-1,1,-1,1,1)$$
$$,(-1,1,-1,-1,1,-1,1,1,-1).$$
Hence  $\mathcal{E}(\Lambda U_{1},\delta_{1})\cong \Bbb
Z_{2}\oplus\Bbb
 Z_{2}\oplus\Bbb Z_{2}\oplus\Bbb
 Z_{2}$. Now let  $(\Lambda U_{2},\delta_{2})$ be the  mcca obtained
from $(\Lambda U_{1},\delta_{1})$  by adding a new generator $x_{4}$
with
 $|x_{4}|=4$ and where the differential is $\delta_{2}(x_{4})=0$, $\delta_{2}(z)=\delta_{1}(z)+x_{4}^{30}$ and $\delta_{2}=\delta_{1}$ on the other generators.
If we write $\xi^{4}:U_{2}^{4}=\Bbb Q\{x_{4}\}\to U_{2}^{4}=\Bbb
Q\{x_{4}\}$ as $\xi^{4}(x_{4})=p_{4}x_{4}$ with $p_{4}\in \Bbb Q$.
Then  we find the same  equations given in (\ref{1}) and the
relation $p_{4}^{30}=p_{3}^{40}$ which have  the following
nontrivial solutions:
$$p_{2}=p_{3}=p_{4}=p_{12}=p_{41}=p_{45}=\pm
1\,\,\,\,\,\,\,\,,\,\,\,\,\,\,\,\,p_{10}=p_{43}=p_{119}=1.$$ Hence
we  get  16 coherent automorphisms of order 2. So
$\mathcal{E}(\Lambda U_{2},\delta_{2})\cong
\underset{2^{4}}{\oplus}\Bbb Z_{2}$ ($2^{4}$ copies of $\Bbb
Z_{2}$)}.

\noindent \textnormal{Next $(\Lambda U_{3},\delta_{3})$ is the  mcca
obtained from $(\Lambda U_{2},\delta_{2})$  by adding a new
generator $x_{5}$ with
 $|x_{5}|=5$ and where the differential is $\delta_{3}(x_{5})=0$, $\delta_{3}(z)=\delta_{2}(z)+x_{5}^{20}$ and $\delta_{3}=\delta_{2}$ on the other generators.
If we write $\xi^{5}:U_{3}^{5}=\Bbb Q\{x_{5}\}\to U_{3}^{5}=\Bbb
Q\{x_{5}\}$ as $\xi^{5}(x_{5})=p_{5}x_{5}$ with $p_{5}\in \Bbb Q$,
Then  we find the same  equations given in (\ref{1}) and the
relations $p_{5}^{20}=p_{4}^{30}$ which have  the following
nontrivial solutions:
$$p_{2}=p_{3}=p_{4}=p_{5}=p_{12}=p_{41}=p_{45}=\pm
1\,\,\,\,\,\,\,\,,\,\,\,\,\,\,\,\,p_{10}=p_{43}=p_{119}=1.$$ So we
get 32 coherent automorphisms of order 2 and $\mathcal{E}(\Lambda
U_{3},\delta_{3})\cong \underset{2^{5}}{\oplus}\Bbb Z_{2}$}.

\noindent \textnormal{Now define the following  mccas}:

\textnormal{ $(\Lambda U_{4},\delta_{4})$ is the  mcca
obtained from $(\Lambda U_{3},\delta_{3})$  by adding a new
generator $x_{6}$ with
 $|x_{6}|=6$ and where the differential is $\delta_{4}(x_{6})=0$, $\delta_{4}(z)=\delta_{3}(z)+x_{6}^{20}$ and $\delta_{4}=\delta_{3}$
 on the other generators.}

$(\Lambda U_{5},\delta_{5})$ is the  mcca
obtained from $(\Lambda U_{4},\delta_{4})$  by adding a new
generator $x_{15}$ with
 $|x_{15}|=15$ and where the differential is $\delta_{5}(x_{15})=0$, $\delta_{5}(z)=\delta_{4}(z)+x_{15}^{8}$ and $\delta_{5}=\delta_{4}$
 on the other generators.\\
 
  $(\Lambda U_{6},\delta_{6})$ is the  mcca obtained from
$(\Lambda U_{5},\delta_{5})$  by adding a new generator $x_{20}$
with
 $|x_{20}|=20$ and where the differential is $\delta_{6}(x_{20})=0$, $\delta_{6}(z)=\delta_{5}(z)+x_{20}^{6}$ and $\delta_{6}=\delta_{5}$
 on the other generators.\\
 
 $(\Lambda U_{7},\delta_{7})$ is the  mcca obtained
from $(\Lambda U_{6},\delta_{6})$  by adding a new generator
$x_{30}$ with
 $|x_{30}|=30$ and where the differential is $\delta_{7}(x_{6})=0$, $\delta_{7}(z)=\delta_{6}(z)+x_{30}^{4}$ and $\delta_{7}=\delta_{6}$
 on the other generators.\\
 
$(\Lambda U_{8},\delta_{8})$ is the  mcca obtained from $(\Lambda
U_{7},\delta_{7})$  by adding a new generator $x_{60}$ with
 $|x_{60}|=60$ and where the differential is $\delta_{7}(x_{60})=0$, $\delta_{8}(z)=\delta_{7}(z)+x_{60}^{2}$ and $\delta_{8}=\delta_{7}$
 on the other generators. 
 
 By the same arguments developed above we get:
$$\mathcal{E}(\Lambda U_{4},\delta_{4})\cong
\underset{2^{6}}{\oplus}\Bbb Z_{2} \,\,,\,\,\mathcal{E}(\Lambda
U_{5},\delta_{5})\cong \underset{2^{7}}{\oplus}\Bbb
Z_{2}\,\,\,,\,\,\,
  \mathcal{E}(\Lambda U_{6},\delta_{6})\cong
\underset{2^{8}}{\oplus}\Bbb Z_{2}$$
$$\mathcal{E}(\Lambda
U_{7},\delta_{7})\cong \underset{2^{9}}{\oplus}\Bbb
Z_{2}\,\,\,\,\,\,,\,\,\,\mathcal{E}(\Lambda U_{8},\delta_{8})\cong
\underset{2^{10}}{\oplus}\Bbb Z_{2}.$$ 
Therefore the groups
$\underset{2^{n}}{\oplus}\Bbb Z_{2}$ are realizable for $n\leq 10$.
\end{remark}

\textnormal{Finally  we end  this work by  asking  the following
question: Is it true that the groups $\underset{2^{n}}{\oplus}\Bbb
 Z_{2}$ are always realizable  for $n\geq 11$?}

\end{document}